\documentclass[12pt]{amsart}
\usepackage{amsmath,amssymb,amsthm}
\usepackage{accents}
\usepackage[top=30truemm,bottom=30truemm,left=25truemm,right=25truemm]{geometry}

\newcommand{\p}{\text{p}}
\newcommand{\LL}{\mathrm{L}}
\newcommand{\V}{\mathrm{V}}
\newcommand{\PP}{\mathbb{P}}
\newcommand{\R}{\mathbb{R}}
\newcommand{\ZF}{\sf{ZF}}
\newcommand{\ZFC}{\sf{ZFC}}
\newcommand{\DC}{\sf{DC}}
\newcommand{\AD}{\sf{AD}}
\newcommand{\fB}{\mathfrak{B}}
\newcommand{\Baire}{\omega^{\omega}}
\newcommand{\PI}{\PP_I}
\newcommand{\bPI}{\widetilde{\PP}_I}
\newcommand{\gn}{\dot{x}_{\text{gen}}}
\newcommand{\ck}{\check}
\newcommand{\ovl}{\overline}
\newcommand{\bface}[1]{\undertilde{\mathbf{#1}}}

\theoremstyle{definition}
\newtheorem{defn}{Definition}[section]
\newtheorem{thm}[defn]{Theorem}
\newtheorem{lem}[defn]{Lemma}
\newtheorem{rem}[defn]{Remark}
\newtheorem{cor}[defn]{Corollary}
\newtheorem{Q}[defn]{Question}

\begin{document}

\keywords{Regularity properties, Determinacy of infinite games, Descriptive set theory.}
\subjclass{03E15, 03E60, 28A05}



\title[Regularity properties, determinacy, and $\infty$-Borel sets]{$I$-regularity, determinacy, and $\infty$-Borel sets of reals}


\author[D.\ Ikegami]{Daisuke Ikegami}
\address[D.\ Ikegami]{College of Engineering, Shibaura Institute of Technology, 3-7-5 Toyosu, Koto-ward, Tokyo, 135-8548 JAPAN}

\email[D.\ Ikegami]{\href{mailto:ikegami@shibaura-it.ac.jp}{ikegami@shibaura-it.ac.jp}}



\thanks{The author would like to thank the Japan Society for the Promotion of Science (JSPS) for its generous support through the grant with JSPS KAKENHI Grant Number 19K03604. He is also grateful to the Sumitomo Foundation for its generous support through Grant for Basic Science Research.}

\begin{abstract}
  We show under $\ZF + \DC + \AD_\R$ that every set of reals is $I$-regular for any $\sigma$-ideal $I$ on the Baire space $\omega^{\omega}$ such that $\mathbb{P}_I$ is proper. This answers the question of Khomskii~\cite[Question~2.6.5]{Yurii_PhD}. We also show that the same conclusion holds under $\ZF + \DC + \AD^+$ if we additionally assume that the set of Borel codes for $I$-positive sets is $\bface{\Delta}^2_1$. If we do not assume $\DC$, the notion of properness becomes obscure as pointed out by Asper\'{o} and Karagila~\cite{DCproper}. Using the notion of strong properness similar to the one introduced by Bagaria and Bosch~\cite{MR2141462}, we show under $\ZF + \DC_{\R}$ without using $\DC$ that every set of reals is $I$-regular for any $\sigma$-ideal $I$ on the Baire space $\Baire$ such that $\mathbb{P}_I$ is strongly proper assuming every set of reals is $\infty$-Borel and there is no $\omega_1$-sequence of distinct reals. In particular, the same conclusion holds in a Solovay model.
\end{abstract}

\maketitle


\section{Introduction}

Regularity properties for sets of reals have been extensively studied since the early 20th century. A set $A$ of reals has some regularity property when $A$ can be approximated by a simple set (such as a Borel set) modulo some small sets. Typical examples of regularity properties are Lebesgue measurability, the Baire property, the perfect set property, and Ramseyness.

Initially motivated by the study of regularity properties, the theory of infinite games has been developed. The work of Banach and Mazur, Davis, and Gale and Stewart shows that the Axiom of Determinacy ($\AD$) implies every set of reals is Lebesgue measurable and every set of reals has the Baire property and the perfect set property. However, it is still open under $\ZF + \DC$ whether $\AD$ implies that every set of reals is Ramsey.

In 1960s, Solovay~\cite{Solovay_70} proved that if the theory $\ZFC + \lq\lq \text{There is an inaccessible cardinal}"$ is consistent, then so is the theory $\ZF + \DC + \lq\lq \text{Every set of reals is Lebesgue measurable}"$. The model of the latter theory he constructed is nowadays called a Solovay model. He has also shown that in a Solovay model, every set of reals has the Baire property and the perfect set property. Later, Mathias proved that in a Solovay model, every set of reals is Ramsey.

After the above results on $\AD$ and Solovay models, many other regularity properties for sets of reals have been investigated. Ikegami~\cite{fa_reg} developed a general framework of regularity properties by introducing the notion of strongly arboreal forcings and assigning a regularity property called $\PP$-measurability to each strongly arboreal forcing $\PP$. Many of the regularity properties are equivalent to $\PP$-measurability for some $\PP$ and he has shown the general equivalence among $\PP$-measurability for $\bface{\Delta}^1_2$-sets of reals, the generic absoluteness for $\bface{\Sigma}^1_3$-statements via $\PP$, and a transcendence property over $\LL[x]$ for all reals $x$ if $\PP$ is proper and simply definable.

The above framework of strongly arboreal forcings can be subsumed using the notion of idealized forcing introduced by Zapletal~\cite{MR2391923}. Starting with a $\sigma$-ideal $I$ on a Polish space, he introduced the forcing $\PI$ consisting of Borel sets not in $I$ ordered by inclusion modulo $I$. Any strongly arboreal forcing $\PP$ is forcing equivalent to some $\PI$ if $\PP$ is proper. Using this wider framework, Khomskii~\cite{Yurii_PhD} introduced $I$-regularity for any $\sigma$-ideal $I$ on the Baire space $\Baire$. The notion of $I$-regularity generalizes $\PP$-measurability and it captures a wider class of regularity properties.

Related to the work on $\AD$ and regularity properties, Khomskii~\cite[Question~2.6.4]{Yurii_PhD} asked the following question:
\begin{Q}[Khomskii~\cite{Yurii_PhD}]\label{q:Yurii}
Does $\AD$ imply every set of reals is $I$-regular for any $\sigma$-ideal on the Baire space $\Baire$ such that the forcing $\PI$ is proper?
\end{Q}

The positive answer to Question~\ref{q:Yurii} would give us that $\AD$ implies every set of reals is Ramsey. So solving Question~\ref{q:Yurii} in a positive way may be difficult. Considering this point, Khomskii~\cite[Question~2.6.5]{Yurii_PhD} asked the following question as well:
\begin{Q}[Khomskii~\cite{Yurii_PhD}]\label{q:Yurii2}
Does $\AD_\R$ imply every set of reals is $I$-regular for any $\sigma$-ideal on the Baire space $\Baire$ such that the forcing $\PI$ is proper?
\end{Q}

In this paper, we will give a positive answer to Question~\ref{q:Yurii2} as follows:
\begin{thm}\label{thm:ADR}
Assume $\ZF + \DC + \AD_\R$. Then for any $\sigma$-ideal $I$ on the Baire space $\Baire$ such that $\PI$ is proper, every set of reals is $I$-regular.
\end{thm}

Using the reflection argument on $\AD^+$, we will show that the assumption of $\AD_\R$ in Theorem~\ref{thm:ADR} can be replaced by $\AD^+$ if the ideal $I$ is simply definable:
\begin{thm}\label{thm:ADplus}
Assume $\ZF + \DC + \AD^+$. Let $I$ be a $\sigma$-ideal on the Baire space $\Baire$ such that $\PI$ is proper and the set $\bPI = \{ c \mid B_c \notin I \}$ is $\undertilde{\mathbf{\Delta}}^2_1$. Then every set of reals is $I$-regular.
\end{thm}

If we do not assume $\DC$, the notion of properness becomes obscure as pointed out by Asper\'{o} and Karagila~\cite{DCproper}. Using the notion of strong properness in Definition~\ref{def:strong properness} similar to the one introduced by Bagaria and Bosch~\cite[Definition~5]{MR2141462}, we will prove the following theorem:
\begin{thm}\label{thm:infty-Borel}
Assume $\ZF + \DC_\R$ and every set of reals is $\infty$-Borel. Suppose also that there is no $\omega_1$-sequence of distinct reals. Let $I$ be a $\sigma$-ideal on the Baire space $\Baire$ such that $\bPI$ is strongly-proper. Then every set of reals is $I$-regular.
\end{thm}

It is not difficult to see that in a Solovay model $\V (\R^*)$, every set of reals is $\infty$-Borel and there is no $\omega_1$-sequence of distinct reals. Therefore, we have the following:
\begin{cor}
In a Solovay model, the following holds: Let $I$ be a $\sigma$-ideal on the Baire space $\Baire$ such that $I$ is strongly proper. Then every set of reals is $I$-regular.
\end{cor}

Section~\ref{sec:basicnotions} is devoted to introducing basic notions and theorems we will use throughout this paper. In Section~\ref{sec:ADR}, we prove Theorem~\ref{thm:ADR}. We prove Theorem~\ref{thm:ADplus} in Section~\ref{sec:ADplus} and Theorem~\ref{thm:infty-Borel} in Section~\ref{sec:infty-Borel}. We end the paper with raising some questions in Section~\ref{sec:Q}.

\section{Basic notions}\label{sec:basicnotions}

From now on, we work in $\ZF+\DC_{\mathbb{R}}$, where $\DC_\R$ states that any relation on the reals with no minimal element has an $\omega$-descending chain.
We assume that readers are familiar with the basics of forcing and descriptive set theory. For basic definitions not given in this paper, see Jech~\cite{Jech} and Moschovakis~\cite{new_Moschovakis}.
When we say \lq \lq reals", we mean elements of the Baire space or of the Cantor space. By $\fB (\Baire)$, we mean the collection of Borel subsets of the Baire space $\Baire$.

In this section, we introduce basic notions and notations we will use throughout this paper.
We start with the central notion of this paper, $I$-regularity for a $\sigma$-ideal $I$ on the Baire space $\Baire$.
\begin{defn}\label{def:I-regularity}
Let $I$ be a $\sigma$-ideal on the Baire space $\Baire$.
\begin{enumerate}
\item A subset $A$ of the Baire space $\Baire$ is called {\it $I$-positive} if $A$ is not in $I$.

\item Let $\PI$ be the collection of Borel sets which are $I$-positive, i.e., $\PI = \{ B \in \fB (\Baire) \mid \text{$B$ is $I$-positive}\}$. For $B, C \in \PI$, $C \le_I B$ if $C \setminus B \in I$.

\item Let $\bPI$ be the collection of Borel codes whose decoded Borel sets are $I$-positive, i.e., $\bPI = \{ c \mid B_c \in \PI\}$. For $c, d \in \bPI$, $c \le_I d$ if $B_c \le_I B_d$.

\item (Khomskii) A subset $A$ of the Baire space $\Baire$ is {\it $I$-regular} if for any $B$ in $\PI$, there is a $C \le_I B$ such that either $C \subseteq A$ or $C \cap A = \emptyset$.
\end{enumerate}
\end{defn}

It is clear that $\PI$ and $\bPI$ are forcing equivalent. We often confuse $\bPI$ with $\PI$ while we use $\bPI$ when we consider the notion of strong properness in Definition~\ref{def:strong properness} (cf.~Remark~\ref{rem:strong properness}).

Many regularity properties for sets of reals can be expressed as $I$-regularity for some $I$. For example, Lebesgue measurability coincides with $I$-regularity when $I$ is the ideal of Lebesgue null sets, and the Baire property is the same as $I$-regularity when $I$ is the ideal of meager sets. If $I$ is the Ramsey null ideal, $I$-regularity is the same as complete Ramseyness. More examples can be found in Khomskii~\cite[Table~2.1]{Yurii_PhD}.

Asper\'{o} and Karagila~\cite{DCproper} modified the definition of hereditary sets $H(\kappa)$ in such a way that it can be defined in $\ZF$ without using the Axiom of Choice while ensuring some basic facts on $H (\kappa)$, and that it is equivalent to the standard definition of $H(\kappa)$ under $\ZFC$. Using this modified definition of $H(\kappa)$, they developed the basic theory of proper forcings under $\ZF + \DC$.
\begin{defn}[Asper\'{o} and Karagila~\cite{DCproper}]\label{def:proper}
\begin{enumerate}
  \item Given an infinite cardinal $\kappa$, let $H(\kappa)$ be the collection of all the sets $x$ such that there is no surjection from the transitive closure of $x$ to $\kappa$.

  \item Assume $\DC$. Let $\PP$ be a poset. We say $\PP$ is {\it proper} if for any sufficiently large cardinal $\kappa$ and every countable elementary substructure $X$ of $(H (\kappa)  , \in)$ with $\PP \in X$, if $p \in \PP \cap X$, then there is a condition $q \le p$ such that $q$ is $(X, \PP)$-generic, i.e., for any predense subset $D$ of $\PP$ in $X$, $D\cap X$ is predense below $q$.
\end{enumerate}
\end{defn}


As pointed by Asper\'{o} and Karagila~\cite{DCproper}, if we do not assume $\DC$, the notion of properness becomes obscure (or every forcing would become proper) given that $\DC$ is equivalent to having countable elementary substrctures of $\V_{\alpha}$ for any infinite ordinal $\alpha$.

In the context of $\ZF+ \DC_{\R}$ without assuming $\DC$, instead of properness, we consider a strengthening of properness similar to the one in Bagaria and Bosch~\cite[Definition~5]{MR2141462}.
\begin{defn}\label{def:strong properness}
Let $\PP$ be a poset. We say $\PP$ is {\it strongly proper} if for any countable transitive model $M$ of a fragment of $\ZF + \DC_\R$ such that $\PP \cap M$, $\le \cap M$, and $\bot \cap M$ are in $M$, if $p \in \PP \cap M$, then there is a condition $q \le p$ such that $q$ is $(M,\PP)$-generic, i.e., if $M \vDash \lq\lq \text{$D$ is a predense subset of $\PP \cap M$}"$, then $D \cap M$ is predense below $q$.
\end{defn}

We make some remarks on Definition~\ref{def:strong properness}.
\begin{rem}\label{rem:strong properness}
\begin{enumerate}
\item Let $\PP$ be a poset consisting of reals. Then if $X$ is a countable elementary substructure of $(\V_{\kappa}, \in)$ with $\PP \in X$ and $M$ is the transitive collapse of $X$, then $p$ is $(X,\PP)$-generic if and only if $p$ is $(M,\PP)$-generic. In particular, if $\PP$ is proper, so is strongly proper. In our context, we consider $\bPI$ for such a $\PP$. Here we use $\bPI$ instead of $\PI$ because $\bPI$ consists of reals while $\PI$ consists of Borel sets of reals.

\item All the typical examples of tree-type forcings satisfying Axiom A are strongly proper.

\item In Bagaria and Bosch~\cite{MR2141462}, they consider strong properness only for projective forcings, i.e., the forcings defined in a projective manner with a real parameter. In this paper, we consider a broader class of forcings to include $I$ which are not projectively defined. 

\end{enumerate}
\end{rem}

We will use the following lemmas on $\PI$ and $\bPI$.
\begin{lem}[Zapletal]\label{lem:generic real}
Let $I$ be a $\sigma$-ideal on the Baire space $\Baire$. Then the forcing $\PI$ adds an element $\gn$ of the Baire space such that if $G$ is $\PI$-generic over $\V$, for every Borel set $B \subseteq \Baire$ in $\V$, $B \in G$ if and only if $\gn^G \in B^{\V[G]}$, where $B^{\V[G]}$ is the decode of a Borel code for $B$ calculated in $\V[G]$. In particular, $\V[G] = \V[\gn^G]$.
\end{lem}

\begin{proof}
The arguments in Zapletal~\cite[Proposition~2.1.2]{MR2391923} can be proceeded in $\ZF + \DC_\R$.
\end{proof}

\begin{defn}\label{def:generic real}
Let $I$ be a $\sigma$-ideal on the Baire space $\Baire$.
\begin{enumerate}
\item Let $\kappa$ be a sufficiently large cardinal and $X$ be a countable elementary substructure $X$ of $(\V_\kappa  , \in)$ with $I \in X$. We say an element $x$ of the Baire space $\Baire$ is an {\it $X$-generic real for $I$} if the set $\{ B \in \PI \cap X \mid x \in B\}$ is a filter on $\PI \cap X$ which meets all the predense subsets of the poset $\PI$ that are elements of $X$.

\item Let $M$ be a transitive model of a fragment of $\ZF + \DC_\R$ such that $\PP \cap M$, $\le \cap M$, and $\bot \cap M$ are in $M$. We say an element $x$ of the Baire space $\Baire$ is an {\it $(M, \bPI)$-generic real} if the set $ \{ c \in \bPI^M \mid x \in B_c \}$ is an $\bPI^M$-generic filter over $M$. 
\end{enumerate}
\end{defn}

\begin{lem}[Zapletal]\label{lem:proper char}
We assume $\DC$. Let $I$ be a $\sigma$-ideal on the Baire space $\Baire$. Then the following are equivalent:
\begin{enumerate}
\item the forcing $\PI$ is proper,

\item for any sufficiently large cardinal $\kappa$ and every countable elementary substructure $X$ of $(\V_{\kappa}  , \in )$ with $I \in X$, if $B \in \PI \cap X$, then the set $C  = \{ x \in B \mid x \text{ is $X$-generic for $I$}\}$ is an $I$-positive Borel set.
\end{enumerate}
\end{lem}

\begin{proof}
See Zapletal~\cite[Proposition~2.2.2]{MR2391923}.
\end{proof}

The following lemma can be proven in a similar way to Lemma~\ref{lem:proper char}.
\begin{lem}\label{lem:strong proper char}
Let $I$ be a $\sigma$-ideal on the Baire space $\Baire$. Then the following are equivalent:
\begin{enumerate}
\item the forcing $\bPI$ is strongly proper,

\item for any countable transitive model $M$ of a fragment of $\ZF + \DC_\R$ such that $\bPI^M = \bPI \cap M$, $\le^M = (\le \cap M)$, and $\bot^M = (\bot \cap M)$, if $b \in \bPI \cap N$, then the set $C  = \{ x \in B_b \mid \text{ $x$ is $(M, \bPI)$-generic }\}$ is an $I$-positive Borel set.
\end{enumerate}
\end{lem}


%


We now introduce the key property for sets of reals in this paper, $\infty$-Borelness.
Infinitary Borel codes ($\infty$-Borel codes) are a transfinite generalization of Borel codes: Let $\mathcal{L}_{\infty, 0} (\{\mathbf{a}_{m, n}\}_{m, n\in \omega})$ be the language allowing arbitrary many well-ordered conjunctions and disjunctions and no quantifiers with atomic sentences $\mathbf{a}_{m, n}$ for each $m, n \in \omega$. The {\it $\infty$-Borel codes} are the sentences in $\mathcal{L}_{\infty, 0} (\{\mathbf{a}_{m, n}\}_{m , n \in \omega})$ belonging to any $\Gamma$ such that
\begin{itemize}
\item the atomic sentence $\mathbf{a}_{m, n}$ is in $\Gamma$ for each $m, n\in \omega$,

\item if $\phi $ is in $\Gamma$, then so is $\neg \phi$, and

\item if $\alpha $ is an ordinal and $\langle \phi_{\beta} \mid \beta < \alpha\rangle $ is a sequence of sentences each of which is in $\Gamma$, then $\bigvee_{\beta < \alpha} \phi_{\beta} $ is also in $\Gamma$.
\end{itemize}
To each $\infty$-Borel code $\phi$, we assign a set of reals $B_{\phi}$ in the same way as decoding Borel codes by induction on the construction of $\phi$:
\begin{itemize}
\item if $\phi = \mathbf{a}_{m, n}$, then $B_{\phi} = \{ x \in \Baire \mid x(m) = n\}$,

\item if $\phi = \neg \psi$, then $B_{\phi} = \Baire \setminus B_{\psi}$, and

\item if $\phi = \bigvee_{\beta < \alpha} \psi_{\beta}$, then $B_{\phi} = \bigcup_{\beta < \alpha} B_{\psi_{\beta}}$.
\end{itemize}

\begin{defn}\label{def:infty-Borel}
Let $A$ be a subset of the Baire space $\Baire$. We say $A$ is {\it $\infty$-Borel} if there is an $\infty$-Borel code $\phi$ such that $A = B_{\phi}$.
\end{defn}

As with Borel codes, one can regard $\infty$-Borel codes as wellfounded trees with atomic sentences $\mathbf{a}_{m, n}$ on terminal nodes and decode them by assigning sets of reals on each node recursively from terminal nodes. (If a node has only one successor, then it means \lq \lq negation" and if a node has more than one successor, then it means \lq \lq disjunction".) The only difference between Borel codes and $\infty$-Borel codes is that trees are on $\omega$ for Borel codes while trees are on ordinals for $\infty$-Borel codes. From this visualization, it is easy to see that the statement \lq \lq $\phi$ is an $\infty$-Borel code" is absolute among transitive models of $\ZF$. Also, it can be easily shown that the statement \lq \lq a real $x$ is in $B_{\phi}$" is absolute among transitive models of $\ZF$.

The following characterization of $\infty$-Borel sets is very useful:
\begin{thm}[Folklore]\label{thm:code equivalence}
Let $A$ be a subset of the Baire space $\Baire$. Then the following are equivalent:
\begin{enumerate}
\item $A$ is $\infty$-Borel,

\item there are a first-order formula $\phi$ and a set $S$ of ordinals such that for each real $x$,
\begin{align*}
x \in A \iff \LL[S,x] \vDash \phi (x).
\end{align*}
\end{enumerate}
\end{thm}

\begin{proof}
See Larson~\cite[Theorem~9.0.4]{ADplus}.
\end{proof}



Using the Axiom of Choice, one can easily show that every set of reals is $\infty$-Borel. However, if we do not assume the Axiom of Choice, the notion of $\infty$-Borelness becomes non-trivial. In fact, we will show that every set of reals is $I$-regular for any $\sigma$-ideal $I$ on the Baire such that $\bPI$ is strongly proper if every set of reals is $\infty$-Borel and there is no $\omega_1$-sequence of distinct reals.

We next define $\AD^+$ and introduce some theorems on $\AD^+$.
\begin{defn}\label{def:ADplus}
\begin{enumerate}
\item Let $\Theta$ be the following ordinal: $\Theta = \sup \ \{ \alpha \mid \text{There is a surjection}$ $\pi \colon \mathbb{R} \to \alpha \}$.

\item For any ordinal $\gamma < \Theta$, we consider the product space $\gamma^{\omega}$ whose basic open sets are of the form $[s] = \{ x \in \gamma^{\omega} \mid s \subseteq x\}$, where $s$ is in $\gamma^{<\omega}$.

\item We say {\it $<$$\Theta$-determinacy} holds if for any $\gamma < \Theta$, any continuous $f \colon \gamma^{\omega} \to \Baire$, and any $A \subseteq \Baire$, the subset $f^{-1}(A)$ of $\gamma^{\omega}$ is determined.

\item The axiom {\it $\AD^+$} states that $\AD$, $\DC_\R$, and $<$$\Theta$-determinacy hold and every set of reals is $\infty$-Borel.
\end{enumerate}
\end{defn}

We will list some theorems on $\AD^+$ we will use in this paper. We say a set $A$ of reals is {\it Suslin} if there are an ordinal $\gamma$ and a tree $T$ on $\omega \times \gamma$ such that $A = \p [T]$. We say a set $A$ of reals is {\it co-Suslin} if the compliment $\Baire \setminus A$ is Suslin.
\begin{thm}[Woodin]\label{fact:ADplus}
Assume $\ZF + \AD^+$. Then the following hold.
\begin{enumerate}
\item Every $\bface{\Sigma}^2_1$ statement has a witness which is a $\bface{\Delta}^2_1$-set of reals.

\item Every $\bface{\Sigma}^2_1$ set of reals is Suslin. In particular, every $\bface{\Delta}^2_1$-set of reals is Suslin and co-Suslin.

\item For any subset $A$ of $(\Baire)^{\omega}$ which is Suslin and co-Suslin, the Gale-Stewart game with reals and the payoff set $A$ is determined.
\end{enumerate}
\end{thm}

\begin{proof}
For 1. and 2., see e.g., Steel and Trang~\cite{ADplusreflection}.
For 3., see e.g., Larson~\cite[Section~13]{ADplus}.
\end{proof}

\section{$I$-regularity and $\AD_\R$}\label{sec:ADR}

In this section, we prove Theorem~\ref{thm:ADR}.

\begin{proof}[Proof of Theorem~\ref{thm:ADR}]
Let $A$ be any set of reals. We will show that $A$ is $I$-regular. Let $B$ be any $I$-positive Borel set. We will find a $C \le_I B$ such that either $C \subseteq A$ or $C \cap A = \emptyset$.

Consider the following game $\mathcal{G}(I;A;B)$ which is essentially the same as a Banach-Mazur game for the Stone space of $\PI$. The game $\mathcal{G}(I;A;B)$ is played by two players, player I and player II. They play elements of $\bPI$ in turn, i.e., player I starts with choosing $c_0 \in \bPI$, then player II responds with $c_1 \in \bPI$, then player I moves with $c_2$ and player II chooses $c_3$ and so on. During the game, they need to keep the following conditions:
\begin{enumerate}
\item $B_{c_0} \le_I B$,

\item for any $n\in \omega$, $c_{n+1} \le_I c_n$ and $c_n$ decides the value $\gn (\ck{n})$ in $\bPI$. (We are using the fact that $\PI$ and $\bPI$ are forcing equivalent.)
\end{enumerate}
After $\omega$ moves, they have produced a $\le_I$-descending sequence $(c_n \mid n \in \omega )$ in $\bPI$. Let $y$ be the element of the Baire space $\Baire$ such that for all natural numbers $n$, $c_n \Vdash \lq\lq \gn (\ck{n}) = \ck{y}(\ck{n})"$. Player I wins if $y$ is in $A$ and player II wins if $y$ is not in $A$.

Since each $c_n$ is a real, by $\AD_\R$, the game $\mathcal{G}(I;A;B)$ is determined. We may assume that Player I has a winning strategy $\sigma$ in the game $\mathcal{G}(I;A;B)$. (The case when Player II has a winning strategy in the game can be dealt with in a similar way.) Let $c_0 = \sigma (\emptyset)$.


Let $\kappa$ be a sufficiently large cardinal. Using $\DC$, one can find a countable elementary substructure $X$ of $(\V_{\kappa} , \in )$ such that $I,A, B, \sigma ,c_0   \in X$. Since $\PI$ is proper and $B_{c_0} \in \PI \cap X$, by Lemma~\ref{lem:proper char}, the set $C  = \{ x \in B_{c_0} \mid x \text{ is $X$-generic for $I$}\}$ is an $I$-positive Borel set. Hence $C \in \PI$. Also since $C \subseteq B_{c_0}$, $C \le_I B_{c_0} \le_I B$.

We will show that $C \subseteq A$. Let $x$ be any element of $C$. We will argue that $x$ is in $A$.

By the definition of $C$, $x$ is $X$-generic for $I$, i.e., the set $G_x = \{ B' \in \PI \cap X \mid x \in B'\}$ is a filter on $\PI \cap X$ which meets all the predense subsets of the poset $\PI$ that are elements of $X$.

We now construct a run $(c_n \in X \mid n \in \omega)$ of the game $\mathcal{G}(I;A;B)$ consistent with the strategy $\sigma$ such that for all $n\in \omega$, $B_{c_n} \in G_x$. The real $c_0$ is already given as $\sigma (\emptyset)$ and $B_{c_0}$ is in $G_x$ because $x \in C \subseteq B_{c_0}$. Suppose that $n$ is any natural number and that the sequence $(c_i \mid i < 2n+1)$ has been obtained, it is in $X$, $B_{c_{2n}} \in G_x$, and it is a partial run of the game $\mathcal{G}(I;A;B)$ consistent with the strategy $\sigma$. Then we will find suitable $c_{2n+1}$ and $c_{2n+2}$. Let $\mathcal{D} = \{ B_{\sigma (c_0, \ldots , c_{2n}, c)} \mid c \le_I c_{2n}\}$. Then it is easy to see that the set $\mathcal{D}$ is dense below $B_{c_{2n}}$. Since the sequence $(c_i \mid i < 2n+1)$ and $\sigma$ are in $X$, the set $\mathcal{D}$ is in $X$ as well. Since $c_{2n} \in G_x$ and $x$ is $X$-generic for $I$, the set $G_x$ meets $\mathcal{D}$. Let $B'$ be in $\mathcal{D} \cap G_x$. Since both $\mathcal{D}$ and $B'$ are in $X$, by the definition of $\mathcal{D}$ and elementarity of $X$, there is a $c_{2n+1}$ in $X$ such that $B' = B_{ \sigma (c_0, \ldots , c_{2n}, c_{2n+1})}$. Let $c_{2n+2} = \sigma (c_0, \ldots , c_{2n}, c_{2n+1})$. Then since $\sigma$ is in $X$, the sequence $(c_0, \ldots , c_{2n+2})$ is in $X$ and $B_{c_{2n+2}} = B' \in G_x$. Therefore, we have obtained the desired $c_{2n+1}$ and $c_{2n+2}$.


Let $y$ be the real associated with the run $(c_n \mid n \in \omega)$ of the game $\mathcal{G}(I;A;B)$, i.e., for all $n \in \omega$, $c_n \Vdash \lq\lq \gn (\ck{n}) = \ck{y} (\ck{n})" $. Since the run $(c_n \mid n \in \omega)$ is consistent with $\sigma$ and $\sigma$ is winning for Player I, by the rule of the game $\mathcal{G}(I;A;B)$, we have $y \in A$. So it is enough to show that $y =x$.

We will verify that $y =x$. Let $n$ be any natural number and $m = y (n)$. We will show that $x(n) = m$. 
Since $c_n \Vdash \lq\lq \gn (\ck{n}) = \ck{y} (\ck{n})"$ and $c_n \in G_x$, by the genericity of $G_x$ over $X$, $y (n) = \gn^{G_x} (n) $, so $\gn^{G_x} (n) = m$. Now 
\begin{align*}
\gn^{G_x}(n) = m	\iff	&	\gn^{G_x} \in \{ z \mid z (n) = m \}\\
		\iff	&	\{ z \mid z (n) = m\} \in G_x\\
		\iff	&	x \in \{ z \mid z (n) = m\}\\
		\iff	& 	x(n) =m
\end{align*}
Hence, $\gn^{G_x}(n) = m = x (n)$. Since $n$ is an arbitrary natural number and $y(n) = m = x(n)$, we have $y = x$.



Since $y$ is in $A$ and $y =x$, the real $x$ is in $A$ as well. Since $x$ was an arbitrary element of $C$, it follows that $C \subseteq A$. Hence we have found a $C \le_I B$ such that either $C \subseteq A$ or $C \cap A = \emptyset$. Since $B$ was an arbitrary $I$-positive Borel set, we have shown that $A$ is $I$-regular. Since $A$ was an arbitrary set of reals, this completes the proof of the theorem.
\end{proof}

\begin{rem}
The assumption \lq\lq $\PI$ is proper" in Theorem~\ref{thm:ADR} is necessary. For example, let $A$ be a $\bface{\Sigma}^1_1$-set of reals which is not Borel and $I = \{ B \in \mathfrak{B}(\Baire) \mid B\cap A \text{ is Borel} \}$. Then $I$ is a non-trivial $\sigma$-ideal on the Baire space $\Baire$ and $A$ is not $I$-regular. This is because for any $B \in \PI$, both $B \cap A$ and $B \setminus A$ are not Borel, in particular neither $B \subseteq A$ nor $B\cap A = \emptyset$ holds. (For the proof of non-properness of this $\PI$, see Zapletal~\cite[Example~2.2.3]{MR2391923}.)
\end{rem}

\section{$I$-regularity and $\AD^+$}\label{sec:ADplus}

In this section, we prove Theorem~\ref{thm:ADplus}.

Notice that the assumption \lq\lq the set $\bPI = \{ c \mid B_c \notin I \}$ is $\undertilde{\mathbf{\Delta}}^2_1$" is harmless because all the typical examples of $I$ satisfy this condition.

\begin{proof}[Proof of Theorem~\ref{thm:ADplus}]
Suppose not. Then there is a set $A$ of reals such that $A$ is not $I$-regular. We will derive a contradiction from this assumption.

First note that the statement \lq\lq There is a set $A$ of reals which is not $I$-regular" is $\bface{\Sigma}^2_1$ indicated as follows:
\begin{align*}
&\text{There is a set $A$ of reals which is not $I$-regular.} \\
\iff &(\exists A \subseteq \Baire )\ (\exists c \in \Baire) \ \big[ c \in \bPI \text{ and }\\
& (\forall d \in \Baire) \ \text{ if $d \in \bPI$ and $d \le_I c$,}
\text{ then $B_d \nsubseteq A$ and $B_d \cap A \neq \emptyset$.} \big]
\end{align*}

Therefore, by item 1.\! of Theorem~\ref{fact:ADplus}, there is a $\bface{\Delta}^2_1$ set $A$ of reals such that $A$ is not $I$-regular. By the construction of $\gn$ in Zapletal~\cite[Proposition~2.1.2]{MR2391923}, the name $\gn$ is simply definable from $\bPI$. Therefore, it follows that for any $B \in \PI$, the payoff set $P \subseteq (\Baire)^{\omega}$ of the game $\mathcal{G}(I;A;B)$ in the proof of Theorem~\ref{thm:ADR} is $\bface{\Delta}^2_1$. Hence by items 2.\! and 3.\! in Theorem~\ref{fact:ADplus}, $\AD^+$ implies that the game $\mathcal{G}(I;A;B)$ is determined for every $B \in \PI$. The proof of Theorem~\ref{thm:ADR} shows that the determinacy of the games $\mathcal{G}(I;A;B)$ for all $B \in \PI$ together with properness of $\PI$ implies that $A$ is $I$-regular. This is a contradiction. Therefore, every set of reals is $I$-regular.
\end{proof}

\section{$I$-regularity and $\infty$-Borel sets}\label{sec:infty-Borel}

In this section, we prove Theorem~\ref{thm:infty-Borel}.
\begin{proof}[Proof of Theorem~\ref{thm:infty-Borel}]
Let $A$ be any set of reals. We will show that $A$ is $I$-regular. Let $B$ be any $I$-positive Borel set. We will find a $C \le_I B$ such that either $C \subseteq A$ or $C \cap A = \emptyset$.

Since every set of reals is $\infty$-Borel, there are first-order formulas $\phi$ and $\psi$ and sets of ordinals $S$ and $T$ such that for all reals $x$ and $c$,
\begin{align}
x \in A \iff \LL[S, x] \vDash \phi [S, x]	\hspace{1cm}	\text{ and }	\hspace{1cm}	c \in \bPI \iff \LL[T, c] \vDash \psi [T, c].
\end{align}

Let $b$ be a Borel code for $B$, i.e., $B_b = B$ and $N = \LL[S,T, b]$. Then by (1) above, the sets $\bPI \cap N$, $\le_{\bPI} \!\!\! \cap N$, and  $\bot_{\bPI} \! \cap N$ are all in $N$. Using this, letting $I^N$ be the $\sigma$-ideal in $N$ generated by the family of sets $\{ B_c^N \mid c \in \bPI \cap N \}$, it is easy to verify that $N$ thinks that $(\bPI \cap N, \le_{\bPI}\!\! \cap N)$ and $(\PP_{I})^N$ are forcing equivalent.

Let $b' \in \bPI \cap N$ be such that $B_{b'} \le_I B_b \, (= B)$ and $N$ thinks either $ B_{b'} \Vdash_{\PP_{I}} \lq\lq \LL[\ck{S}, \gn] \vDash \phi [\ck{S}, \gn]"$ or $B_{b'} \Vdash_{\PP_{I}}\lq\lq \LL[\ck{S}, \gn] \vDash \neg \phi [\ck{S}, \gn]"$. We may assume the former case $N \vDash B_{b'} \Vdash_{\PP_{I}} \lq\lq \LL[\ck{S}, \gn] \vDash \phi [\ck{S}, \gn]"$. (The latter case $N \vDash B_{b'} \Vdash_{\PP_{I}}\lq\lq \LL[\ck{S}, \gn] \vDash \neg \phi [\ck{S}, \gn]"$ can be dealt with in a similar way to the arguments below.)

Let $\kappa$ be a sufficiently large cardinal. Since $N$ is an inner model of $\ZFC$, we can find an elementary substructure $X$ of $(\V_{\kappa}, \in)^N$ such that $\R^N , \mathcal{P}(\R)^N \subseteq X$, $S,T, b, b' \in X$, and that in $N$, the cardinality of $X$ is the same as that of $\R^N \cup \mathcal{P}(\R)^N$. Since there is no $\omega_1$-sequence of distinct reals and $N$ is an inner model of $\ZFC$, the sets $\mathbb{R}^N$ and $\mathcal{P}(\R)^N$ are both countable in $\V$. Therefore, $X$ is also countable in $\V$.

Let $\pi \colon X \to M$ be the Mostowski collapsing map of $X$. From now on, we write $\ovl{a}$ for $\pi (a)$ for any element $a$ of $X$. Since $\mathbb{R}^N$ and $\mathcal{P}(\R)^N$ are subsets of $X$, it follows that $\bPI \cap N = \ovl{\bPI \cap N} = \bPI \cap M$ is in $M$. Similarly $\le_{\bPI} \cap M$ and $\bot_{\bPI} \cap M$ are in $M$. Since $M$ is a countable transitive model of a fragment of $\ZFC$ and $b' = \ovl{b'} \in M$, by strong properness of $\bPI$ and Lemma~\ref{lem:strong proper char}, the set $C = \{ x \in B_{b'} \mid \text{ $x$ is $(M, \bPI)$-generic}\}$ is an $I$-positive Borel set. Hence $C \in \PI$ and $C \le_I B_{b'} \le_I B$.

We will show that $C \subseteq A$. Let $x$ be any element of $C$. We will argue that $x$ is in $A$.

We first claim that $x$ is $(N, \bPI)$-generic.
Since $x$ is $(M, \bPI)$-generic, the set $G'_x = \{ c \in \bPI \cap M \mid x \in B_c\}$ forms a filter on $\bPI^M$ which meets all the predense open subsets of the poset $\bPI^M = \bPI \cap M$ that are elements of $M$. Since $\bPI \cap M = \bPI \cap N$ and $\mathcal{P}(\R)^M = \mathcal{P}(\R)^N$, the filter $G'_x$ meets all the predense subsets of the poset $\bPI^N = \bPI \cap N$ that are elements of $N$ as well. Hence $x$ is $(N, \bPI)$-generic.

Since $N$ thinks that $(\bPI \cap N, \le_{\bPI}\!\! \cap N)$ and $(\PP_{I})^N$ are forcing equivalent, letting $G_x = \{ B \in (\PP_{I})^N \mid x \in B\}$, $G_x$ is a $(\PP_{I})^N$-generic filter over $N$.
Since $N$ thinks $ B_{b'} \Vdash_{\PP_{I}} \lq\lq \LL[\ck{S}, \gn] \vDash \phi [\ck{S}, \gn]"$, $N[G_x]$ thinks that $\LL[S,\gn^{G_x} ] \vDash \phi [S, \gn^{G_x}]$ and hence $\LL[S,\gn^{G_x} ] \vDash \phi [S, \gn^{G_x}]$ in $\V$ as well. By (1) above, this gives us that the real $\gn^{G_x}$ is in $A$.

We will finish arguing that $x$ is in $A$ by verifying $x = \gn^{G_x}$. Applying Lemma~\ref{lem:generic real} in $N$ to $({\PP}_{I})^N$ and $G_x$, for any Borel set $B''$ in $N$, $x \in (B'')^{N[G_x]} \iff B'' \in G_x \iff \gn^{G_x} \in (B'')^{N[G_x]}$. In particular, $\gn^{G_x} =x$. Therefore, $x = \gn^{G_x} \in A$.

Since $x$ was an arbitrary element of $C$, it follows that $C \subseteq A$. So we have found a $C \le_I B$ such that either $C \subseteq A$ or $C \cap A = \emptyset$. Since $B$ was an arbitrary $I$-positive Borel set, we have shown that $A$ is $I$-regular. Since $A$ was an arbitrary set of reals, this completes the proof of the theorem.
\end{proof}

\section{Questions}\label{sec:Q}

We end this paper with some questions.
\begin{Q}\label{q1}
Does $\AD$ imply every set of reals is $I$-regular for any $\sigma$-ideal on the Baire space $\Baire$ such that $\PI$ is strongly proper?
\end{Q}

The positive answer to Question~\ref{q1} would give us that $\AD$ implies that every set of reals is Ramsey. If $\AD$ implies that every set of reals is $\infty$-Borel, then by Theorem~\ref{thm:infty-Borel}, we would get the positive answer to Question~\ref{q1}.

\begin{Q}\label{q2}
Assume $\DC$ and every set of reals is $\infty$-Borel. Suppose also that there is no $\omega_1$-sequence of distinct reals. Let $I$ be a $\sigma$-ideal on the Baire space $\Baire$ such that $\PI$ is proper. Then must every set of reals be $I$-regular?
\end{Q}

In Khomskii~\cite[Proposition~2.2.8]{Yurii_PhD}, it is claimed that in a Solovay model, every set of reals is $I$-regular for any $\sigma$-ideal $I$ on the Baire space $\Baire$ such that $\PI$ is proper. Unfortunately, there is a gap in the proof there.

\begin{Q}\label{q3}
Does a Solovay model satisfy the following: Let $I$ be a $\sigma$-ideal on the Baire space $\Baire$ such that $\PI$ is proper. Then every set of reals is $I$-regular.
\end{Q}

The positive answer to Question~\ref{q2} would give a positive answer to Question~\ref{q3}.

\newpage

\bibliographystyle{plain}
\bibliography{myreference}

\end{document}